\def\temp{1.34}%
\let\tempp=\relax
\expandafter\ifx\csname psboxversion\endcsname\relax
  \message{PSBOX(\temp) loading}%
\else
    \ifdim\temp cm>\psboxversion cm
      \message{PSBOX(\temp) loading}%
    \else
      \message{PSBOX(\psboxversion) is already loaded: I won't load
        PSBOX(\temp)!}%
      \let\temp=\psboxversion
      \let\tempp= 
    \fi
\fi
\tempp
\let\psboxversion=\temp
\catcode`\@=11
%
%
\def\psfortextures{
\def\PSspeci@l##1##2{%
\special{illustration ##1\space scaled ##2}%
}}%
\def\psfordvitops{
\def\PSspeci@l##1##2{%
\special{dvitops: import ##1\space \the\drawingwd \the\drawinght}%
}}%
\def\psfordvips{
\def\PSspeci@l##1##2{%
\d@my=0.1bp \d@mx=\drawingwd \divide\d@mx by\d@my
\includegraphics{##1\space}}}%
\def\psforoztex{
\def\PSspeci@l##1##2{%
\special{##1 \space
      ##2 1000 div dup scale
      \number-\psllx\space \number-\pslly\space translate
}}}%
\def\psfordvitps{
\def\psdimt@n@sp##1{\d@mx=##1\relax\edef\psn@sp{\number\d@mx}}
\def\PSspeci@l##1##2{%
\special{dvitps: Include0 "psfig.psr"}
\psdimt@n@sp{\drawingwd}
\special{dvitps: Literal "\psn@sp\space"}
\psdimt@n@sp{\drawinght}
\special{dvitps: Literal "\psn@sp\space"}
\psdimt@n@sp{\psllx bp}
\special{dvitps: Literal "\psn@sp\space"}
\psdimt@n@sp{\pslly bp}
\special{dvitps: Literal "\psn@sp\space"}
\psdimt@n@sp{\psurx bp}
\special{dvitps: Literal "\psn@sp\space"}
\psdimt@n@sp{\psury bp}
\special{dvitps: Literal "\psn@sp\space startTexFig\space"}
\special{dvitps: Include1 "##1"}
\special{dvitps: Literal "endTexFig\space"}
}}%
\def\psfordvialw{
\def\PSspeci@l##1##2{
\special{language "PostScript",
position = "bottom left",
literal "  \psllx\space \pslly\space translate
  ##2 1000 div dup scale
  -\psllx\space -\pslly\space translate",
include "##1"}
}}%
\def\psforptips{
\def\PSspeci@l##1##2{{
\d@mx=\psurx bp
\advance \d@mx by -\psllx bp
\divide \d@mx by 1000\multiply\d@mx by \xscale
\incm{\d@mx}
\let\tmpx\dimincm
\d@my=\psury bp
\advance \d@my by -\pslly bp
\divide \d@my by 1000\multiply\d@my by \xscale
\incm{\d@my}
\let\tmpy\dimincm
\d@mx=-\psllx bp
\divide \d@mx by 1000\multiply\d@mx by \xscale
\d@my=-\pslly bp
\divide \d@my by 1000\multiply\d@my by \xscale
\at(\d@mx;\d@my){\special{ps:##1 x=\tmpx, y=\tmpy}}
}}}%
\def\psonlyboxes{
\def\PSspeci@l##1##2{%
\at(0cm;0cm){\boxit{\vbox to\drawinght
  {\vss\hbox to\drawingwd{\at(0cm;0cm){\hbox{({\tt##1})}}\hss}}}}
}}%
\def\psloc@lerr#1{%
\let\savedPSspeci@l=\PSspeci@l%
\def\PSspeci@l##1##2{%
\at(0cm;0cm){\boxit{\vbox to\drawinght
  {\vss\hbox to\drawingwd{\at(0cm;0cm){\hbox{({\tt##1}) #1}}\hss}}}}
\let\PSspeci@l=\savedPSspeci@l
}}%
%
%
\newread\pst@mpin
\newdimen\drawinght\newdimen\drawingwd
\newdimen\psxoffset\newdimen\psyoffset
\newbox\drawingBox
\newcount\xscale \newcount\yscale \newdimen\pscm\pscm=1cm
\newdimen\d@mx \newdimen\d@my
\newdimen\pswdincr \newdimen\pshtincr
\let\ps@nnotation=\relax
{\catcode`\|=0 |catcode`|\=12 |catcode`|
|catcode`#=12 |catcode`*=14
|xdef|backslashother{\}*
|xdef|percentother{
|xdef|tildeother{~}*
|xdef|sharpother{#}*
}%
\def\R@moveMeaningHeader#1:->{}%
\def\uncatcode#1{%
\edef#1{\expandafter\R@moveMeaningHeader\meaning#1}}%
\def\execute#1{#1}
\def\psm@keother#1{\catcode`#112\relax}
\def\executeinspecs#1{%
\execute{\begingroup\let\do\psm@keother\dospecials\catcode`\^^M=9#1\endgroup}}%
\def\@mpty{}%
\def\matchexpin#1#2{
  \fi%
  \edef\tmpb{{#2}}%
  \expandafter\makem@tchtmp\tmpb%
  \edef\tmpa{#1}\edef\tmpb{#2}%
  \expandafter\expandafter\expandafter\m@tchtmp\expandafter\tmpa\tmpb\endm@tch%
  \if\match%
}%
\def\matchin#1#2{%
  \fi%
  \makem@tchtmp{#2}%
  \m@tchtmp#1#2\endm@tch%
  \if\match%
}%
\def\makem@tchtmp#1{\def\m@tchtmp##1#1##2\endm@tch{%
  \def\tmpa{##1}\def\tmpb{##2}\let\m@tchtmp=\relax%
  \ifx\tmpb\@mpty\def\match{YN}%
  \else\def\match{YY}\fi%
}}%
\def\incm#1{{\psxoffset=1cm\d@my=#1
 \d@mx=\d@my
  \divide\d@mx by \psxoffset
  \xdef\dimincm{\number\d@mx.}
  \advance\d@my by -\number\d@mx cm
  \multiply\d@my by 100
 \d@mx=\d@my
  \divide\d@mx by \psxoffset
  \edef\dimincm{\dimincm\number\d@mx}
  \advance\d@my by -\number\d@mx cm
  \multiply\d@my by 100
 \d@mx=\d@my
  \divide\d@mx by \psxoffset
  \xdef\dimincm{\dimincm\number\d@mx}
}}%
%
\newif\ifNotB@undingBox
\newhelp\PShelp{Proceed: you'll have a 5cm square blank box instead of
your graphics (Jean Orloff).}%
\def\s@tsize#1 #2 #3 #4\@ndsize{
  \def\psllx{#1}\def\pslly{#2}%
  \def\psurx{#3}\def\psury{#4}
  \ifx\psurx\@mpty\NotB@undingBoxtrue
  \else
    \drawinght=#4bp\advance\drawinght by-#2bp
    \drawingwd=#3bp\advance\drawingwd by-#1bp
  \fi
  }%
\def\sc@nBBline#1:#2\@ndBBline{\edef\p@rameter{#1}\edef\v@lue{#2}}%
\def\g@bblefirstblank#1#2:{\ifx#1 \else#1\fi#2}%
{\catcode`\%=12
\xdef\B@undingBox{
\def\ReadPSize#1{
 \readfilename#1\relax
 \let\PSfilename=\lastreadfilename
 \openin\pst@mpin=#1\relax
 \ifeof\pst@mpin \errhelp=\PShelp
   \errmessage{I haven't found your postscript file (\PSfilename)}%
   \psloc@lerr{was not found}%
   \s@tsize 0 0 142 142\@ndsize
   \closein\pst@mpin
 \else
   \if\matchexpin{\GlobalInputList}{, \lastreadfilename}%
   \else\xdef\GlobalInputList{\GlobalInputList, \lastreadfilename}%
     \immediate\write\psbj@inaux{\lastreadfilename,}%
   \fi%
   \loop
     \executeinspecs{\catcode`\ =10\global\read\pst@mpin to\n@xtline}%
     \ifeof\pst@mpin
       \errhelp=\PShelp
       \errmessage{(\PSfilename) is not an Encapsulated PostScript File:
           I could not find any \B@undingBox: line.}%
       \edef\v@lue{0 0 142 142:}%
       \psloc@lerr{is not an EPSFile}%
       \NotB@undingBoxfalse
     \else
       \expandafter\sc@nBBline\n@xtline:\@ndBBline
       \ifx\p@rameter\B@undingBox\NotB@undingBoxfalse
         \edef\t@mp{%
           \expandafter\g@bblefirstblank\v@lue\space\space\space}%
         \expandafter\s@tsize\t@mp\@ndsize
       \else\NotB@undingBoxtrue
       \fi
     \fi
   \ifNotB@undingBox\repeat
   \closein\pst@mpin
 \fi
\message{#1}%
}%
%
%
\def\psboxto(#1;#2)#3{\vbox{%
   \ReadPSize{#3}%
   \advance\pswdincr by \drawingwd
   \advance\pshtincr by \drawinght
   \divide\pswdincr by 1000
   \divide\pshtincr by 1000
   \d@mx=#1
   \ifdim\d@mx=0pt\xscale=1000
         \else \xscale=\d@mx \divide \xscale by \pswdincr\fi
   \d@my=#2
   \ifdim\d@my=0pt\yscale=1000
         \else \yscale=\d@my \divide \yscale by \pshtincr\fi
   \ifnum\yscale=1000
         \else\ifnum\xscale=1000\xscale=\yscale
                    \else\ifnum\yscale<\xscale\xscale=\yscale\fi
              \fi
   \fi
   \divide\drawingwd by1000 \multiply\drawingwd by\xscale
   \divide\drawinght by1000 \multiply\drawinght by\xscale
   \divide\psxoffset by1000 \multiply\psxoffset by\xscale
   \divide\psyoffset by1000 \multiply\psyoffset by\xscale
   \global\divide\pscm by 1000
   \global\multiply\pscm by\xscale
   \multiply\pswdincr by\xscale \multiply\pshtincr by\xscale
   \ifdim\d@mx=0pt\d@mx=\pswdincr\fi
   \ifdim\d@my=0pt\d@my=\pshtincr\fi
   \message{scaled \the\xscale}%
 \hbox to\d@mx{\hss\vbox to\d@my{\vss
   \global\setbox\drawingBox=\hbox to 0pt{\kern\psxoffset\vbox to 0pt{%
      \kern-\psyoffset
      \PSspeci@l{\PSfilename}{\the\xscale}%
      \vss}\hss\ps@nnotation}%
   \global\wd\drawingBox=\the\pswdincr
   \global\ht\drawingBox=\the\pshtincr
   \global\drawingwd=\pswdincr
   \global\drawinght=\pshtincr
   \baselineskip=0pt
   \copy\drawingBox
 \vss}\hss}%
  \global\psxoffset=0pt
  \global\psyoffset=0pt
  \global\pswdincr=0pt
  \global\pshtincr=0pt 
  \global\pscm=1cm 
}}%
%
%
\def\psboxscaled#1#2{\vbox{%
  \ReadPSize{#2}%
  \xscale=#1
  \message{scaled \the\xscale}%
  \divide\pswdincr by 1000 \multiply\pswdincr by \xscale
  \divide\pshtincr by 1000 \multiply\pshtincr by \xscale
  \divide\psxoffset by1000 \multiply\psxoffset by\xscale
  \divide\psyoffset by1000 \multiply\psyoffset by\xscale
  \divide\drawingwd by1000 \multiply\drawingwd by\xscale
  \divide\drawinght by1000 \multiply\drawinght by\xscale
  \global\divide\pscm by 1000
  \global\multiply\pscm by\xscale
  \global\setbox\drawingBox=\hbox to 0pt{\kern\psxoffset\vbox to 0pt{%
     \kern-\psyoffset
     \PSspeci@l{\PSfilename}{\the\xscale}%
     \vss}\hss\ps@nnotation}%
  \advance\pswdincr by \drawingwd
  \advance\pshtincr by \drawinght
  \global\wd\drawingBox=\the\pswdincr
  \global\ht\drawingBox=\the\pshtincr
  \global\drawingwd=\pswdincr
  \global\drawinght=\pshtincr
  \baselineskip=0pt
  \copy\drawingBox
  \global\psxoffset=0pt
  \global\psyoffset=0pt
  \global\pswdincr=0pt
  \global\pshtincr=0pt 
  \global\pscm=1cm
}}%
%
\def\psbox#1{\psboxscaled{1000}{#1}}%
\newif\ifn@teof\n@teoftrue
\newif\ifc@ntrolline
\newif\ifmatch
\newread\j@insplitin
\newwrite\j@insplitout
\newwrite\psbj@inaux
\immediate\openout\psbj@inaux=psbjoin.aux
\immediate\write\psbj@inaux{\string\joinfiles}%
\immediate\write\psbj@inaux{\jobname,}%
%
%
\def\toother#1{\ifcat\relax#1\else\expandafter%
  \toother@ux\meaning#1\endtoother@ux\fi}%
\def\toother@ux#1 #2#3\endtoother@ux{\def\tmp{#3}%
  \ifx\tmp\@mpty\def\tmp{#2}\let\next=\relax%
  \else\def\next{\toother@ux#2#3\endtoother@ux}\fi%
\next}%
%
%
\let\readfilenamehook=\relax
\def\re@d{\expandafter\re@daux}
\def\re@daux{\futurelet\nextchar\stopre@dtest}%
\def\re@dnext{\xdef\lastreadfilename{\lastreadfilename\nextchar}%
  \afterassignment\re@d\let\nextchar}%
\def\stopre@d{\egroup\readfilenamehook}%
\def\stopre@dtest{%
  \ifcat\nextchar\relax\let\nextread\stopre@d
  \else
    \ifcat\nextchar\space\def\nextread{%
      \afterassignment\stopre@d\chardef\nextchar=`}%
    \else\let\nextread=\re@dnext
      \toother\nextchar
      \edef\nextchar{\tmp}%
    \fi
  \fi\nextread}%
\def\readfilename{\bgroup%
  \let\\=\backslashother \let\%=\percentother \let\~=\tildeother
  \let\#=\sharpother \xdef\lastreadfilename{}%
  \re@d}%
%
%
\xdef\GlobalInputList{\jobname}%
\def\psnewinput{%
  \def\readfilenamehook{
    \if\matchexpin{\GlobalInputList}{, \lastreadfilename}%
    \else\xdef\GlobalInputList{\GlobalInputList, \lastreadfilename}%
      \immediate\write\psbj@inaux{\lastreadfilename,}%
    \fi%
    \ps@ldinput\lastreadfilename\relax%
    \let\readfilenamehook=\relax%
  }\readfilename%
}%
\expandafter\ifx\csname @@input\endcsname\relax    
  \immediate\let\ps@ldinput=\input\def\input{\psnewinput}%
\else
  \immediate\let\ps@ldinput=\@@input
  \def\@@input{\psnewinput}%
\fi%
\def\nowarnopenout{%
 \def\warnopenout##1##2{%
   \readfilename##2\relax
   \message{\lastreadfilename}%
   \immediate\openout##1=\lastreadfilename\relax}}%
\def\warnopenout#1#2{%
 \readfilename#2\relax
 \def\t@mp{TrashMe,psbjoin.aux,psbjoint.tex,}\uncatcode\t@mp
 \if\matchexpin{\t@mp}{\lastreadfilename,}%
 \else
   \immediate\openin\pst@mpin=\lastreadfilename\relax
   \ifeof\pst@mpin
     \else
     \errhelp{If the content of this file is so precious to you, abort (ie
press x or e) and rename it before retrying.}%
     \errmessage{I'm just about to replace your file named \lastreadfilename}%
   \fi
   \immediate\closein\pst@mpin
 \fi
 \message{\lastreadfilename}%
 \immediate\openout#1=\lastreadfilename\relax}%
{\catcode`\%=12\catcode`\*=14
\gdef\splitfile#1{*
 \readfilename#1\relax
 \immediate\openin\j@insplitin=\lastreadfilename\relax
 \ifeof\j@insplitin
   \message{! I couldn't find and split \lastreadfilename!}*
 \else
   \immediate\openout\j@insplitout=TrashMe
   \message{< Splitting \lastreadfilename\space into}*
   \loop
     \ifeof\j@insplitin
       \immediate\closein\j@insplitin\n@teoffalse
     \else
       \n@teoftrue
       \executeinspecs{\global\read\j@insplitin to\spl@tinline\expandafter
         \ch@ckbeginnewfile\spl@tinline
       \ifc@ntrolline
       \else
         \toks0=\expandafter{\spl@tinline}*
         \immediate\write\j@insplitout{\the\toks0}*
       \fi
     \fi
   \ifn@teof\repeat
   \immediate\closeout\j@insplitout
 \fi\message{>}*
}*
\gdef\ch@ckbeginnewfile#1
 \def\t@mp{#1}*
 \ifx\@mpty\t@mp
   \def\t@mp{#3}*
   \ifx\@mpty\t@mp
     \global\c@ntrollinefalse
   \else
     \immediate\closeout\j@insplitout
     \warnopenout\j@insplitout{#2}*
     \global\c@ntrollinetrue
   \fi
 \else
   \global\c@ntrollinefalse
 \fi}*
\gdef\joinfiles#1\into#2{*
 \message{< Joining following files into}*
 \warnopenout\j@insplitout{#2}*
 \message{:}*
 {*
 \edef\w@##1{\immediate\write\j@insplitout{##1}}*
\w@{
\w@{
\w@{
\w@{
\w@{
\w@{
\w@{
\w@{
\w@{
\w@{
\w@{\string\input\space psbox.tex}*
\w@{\string\splitfile{\string\jobname}}*
\w@{\string\let\string\autojoin=\string\relax}*
}*
 \expandafter\tre@tfilelist#1, \endtre@t
 \immediate\closeout\j@insplitout
 \message{>}*
}*
\gdef\tre@tfilelist#1, #2\endtre@t{*
 \readfilename#1\relax
 \ifx\@mpty\lastreadfilename
 \else
   \immediate\openin\j@insplitin=\lastreadfilename\relax
   \ifeof\j@insplitin
     \errmessage{I couldn't find file \lastreadfilename}*
   \else
     \message{\lastreadfilename}*
     \immediate\write\j@insplitout{
     \executeinspecs{\global\read\j@insplitin to\oldj@ininline}*
     \loop
       \ifeof\j@insplitin\immediate\closein\j@insplitin\n@teoffalse
       \else\n@teoftrue
         \executeinspecs{\global\read\j@insplitin to\j@ininline}*
         \toks0=\expandafter{\oldj@ininline}*
         \let\oldj@ininline=\j@ininline
         \immediate\write\j@insplitout{\the\toks0}*
       \fi
     \ifn@teof
     \repeat
   \immediate\closein\j@insplitin
   \fi
   \tre@tfilelist#2, \endtre@t
 \fi}*
}%
\def\autojoin{%
 \immediate\write\psbj@inaux{\string\into{psbjoint.tex}}%
 \immediate\closeout\psbj@inaux
 \expandafter\joinfiles\GlobalInputList\into{psbjoint.tex}%
}%
%
%
%
\def\centinsert#1{\midinsert\line{\hss#1\hss}\endinsert}%
\def\psannotate#1#2{\vbox{%
  \def\ps@nnotation{#2\global\let\ps@nnotation=\relax}#1}}%
\def\pscaption#1#2{\vbox{%
   \setbox\drawingBox=#1
   \copy\drawingBox
   \vskip\baselineskip
   \vbox{\hsize=\wd\drawingBox\setbox0=\hbox{#2}%
     \ifdim\wd0>\hsize
       \noindent\unhbox0\tolerance=5000
    \else\centerline{\box0}%
    \fi
}}}%
%
\def\at(#1;#2)#3{\setbox0=\hbox{#3}\ht0=0pt\dp0=0pt
  \rlap{\kern#1\vbox to0pt{\kern-#2\box0\vss}}}%
%
\newdimen\gridht \newdimen\gridwd
\def\gridfill(#1;#2){%
  \setbox0=\hbox to 1\pscm
  {\vrule height1\pscm width.4pt\leaders\hrule\hfill}%
  \gridht=#1
  \divide\gridht by \ht0
  \multiply\gridht by \ht0
  \gridwd=#2
  \divide\gridwd by \wd0
  \multiply\gridwd by \wd0
  \advance \gridwd by \wd0
  \vbox to \gridht{\leaders\hbox to\gridwd{\leaders\box0\hfill}\vfill}}%
%
\def\fillinggrid{\at(0cm;0cm){\vbox{%
  \gridfill(\drawinght;\drawingwd)}}}%
%
%
\def\textleftof#1:{%
  \setbox1=#1
  \setbox0=\vbox\bgroup
    \advance\hsize by -\wd1 \advance\hsize by -2em}%
\def\textrightof#1:{%
  \setbox0=#1
  \setbox1=\vbox\bgroup
    \advance\hsize by -\wd0 \advance\hsize by -2em}%
\def\endtext{%
  \egroup
  \hbox to \hsize{\valign{\vfil##\vfil\cr%
\box0\cr%
\noalign{\hss}\box1\cr}}}%
%
\def\frameit#1#2#3{\hbox{\vrule width#1\vbox{%
  \hrule height#1\vskip#2\hbox{\hskip#2\vbox{#3}\hskip#2}%
        \vskip#2\hrule height#1}\vrule width#1}}%
\def\boxit#1{\frameit{0.4pt}{0pt}{#1}}%
\catcode`\@=12 
%
 \psfordvips   
      
\magnification=\magstep1
\parskip 9pt
\font\rmone= cmr10 scaled \magstep1
\font\ninerm=cmr9
\font\sevenrm=cmr7
\font\Bbb=msbm10                                                               
\textfont12=\Bbb                                                               
\let\mcd=\mathchardef
\mcd\Re="7C52                                                                  
\mcd\Ze="7C5A 
\mcd\Pe="7C50                                                                  
\mcd\Ee="7C45
\def\scr{\scriptstyle}
\def\dis{\displaystyle}
\def\text{\textstyle}
\def\un{\rmone 1}  
\def\ti{\tilde}  
\def\la{\lambda}                                                          
\def\O{{\cal O}}

\centerline{\rmone Hard rods: statistics of parking configurations}
\bigskip\bigskip
\centerline{\bf Fran\c cois Dunlop, Thierry Huillet}
\bigskip
{\ninerm
\centerline{Laboratoire de Physique Th\'eorique et Mod\'elisation - 
CNRS UMR 8089}
\centerline{Universit{\'e} de Cergy-Pontoise, 95031 Cergy-Pontoise, France}
\centerline{\sevenrm dunlop@ptm.u-cergy.fr, huillet@ptm.u-cergy.fr}
}
\bigskip\bigskip
\noindent{\bf Abstract.} We compute the correlation function in
the equilibrium version of R\'enyi's {\sl parking problem}.
The correlation length is found to diverge as  
$2^{-1}\pi^{-2}(1-\rho)^{-2}$ when $\rho\nearrow1$ (maximum density) and as
$\pi^{-2}(2\rho-1)^{-2}$ when $\rho\searrow1/2$ (minimum density).

\noindent {\ninerm PACS codes:  05.20.Jj, 05.70.Ce, 02.50.Cw} 

\noindent {\ninerm KEYWORDS: Hard rods, Parking, Correlation length}

\bigskip\noindent{\bf 1. Introduction}
\smallskip\noindent
A parking configuration of hard rods (or cars) on a line or a circle
is a configuration where the largest gap between successive rods
is less than the length of a rod: the parking is full.
The statistical properties of parking configurations obtained by
random sequential parking (or adsorption)
were first studied by R\'enyi [1]. 
The present paper deals with the corresponding equilibrium model.

The correlation function of a gas of hard rods in one dimension was computed by
Frenkel [2]. A general one-dimensional fluid with a nearest neighbor 
interaction, strongly repulsive at short distance and decaying rapidly
at large distance, was then solved by G\"ursey [3] in a grand canonical 
ensemble. Salsburg, Zwanzig and Kirkwood [4] derived a similar solution in the 
canonical ensemble, suitable for testing the then newly invented 
Kirkwood-Salsburg equations. 

We give a solution which is a little simpler, in a canonical isobaric
ensemble, and give precise asymptotic forms for the diverging correlation
length near maximum density (pressure going to $+\infty$) and near
minimum density  (pressure going to $-\infty$).

The intervals between hard rods can be mapped onto the gradient variables
in a one-dimensional interface model.
One could prove mathematically the equivalence of ensembles for the parking
model, including uniqueness of correlation functions, using the local central 
limit theorem like in the proof of the Wulff shape for one-dimensional 
interfaces [5].

As a by-product, we give the probability for $N$ points distributed uniformly 
and independently in $(0,L)$ to be in a parking configuration, as a function of
$\rho=N/L$, in the thermodynamic limit.
\vfill\eject

\bigskip\noindent{\bf 2. The model}
\smallskip\noindent
A parking configuration of $N$ hard rods (or arcs) of length one on a circle
of length $L$ is specified by $N$ positions $X_1,\dots X_N\in[0,L)=\Re/L\Ze$ 
with
$$
\eqalign{
|X_i-X_j|\ge1\quad\forall i\ne j\cr
\min_{j\ne i}(X_i-X_j)_+\le2\quad\forall i\cr
\min_{j\ne i}(X_j-X_i)_+\le2\quad\forall i
}\eqno{(2.1)}
$$
The Lebesgue measure on the set of parking configurations, a Borel subset
of $\Re^N$, normalised by the total measure of this set, defines a probability
measure $\Pe_{N,L}(\cdot)$. This probability measure inherits the rotation
invariance of the Haar measure on the circle of length $L$.
The canonical partition function $Z_{N,L}$ is defined as
$$
Z_{N,L}={1\over N!}\int_{\rm Parking}dx_1\dots dx_N\eqno{(2.2)}
$$
where the range of integration ``Parking'' is defined by (2.1).
The free energy is $F({N,L})=-\ln Z_{N,L}$. Temperature
plays no role and is omitted. 

For $x,y\in[0,L)$ we define
$$
\eqalign{
\rho^{(1)}_{N,L}(x)=&\sum_{i=1}^N{\Pe_{N,L}(\,X_i\in(x,x+dx)\,)\over dx}\cr
\rho^{(2)}_{N,L}(x,y)=&\sum_{i\ne j=1}^N
{\Pe_{N,L}(\,X_i\in(x,x+dx),X_j\in(y,y+dy)\,)\over dxdy}
}\eqno{(2.3)}
$$
By rotation invariance $\rho^{(1)}_{N,L}(x)$ is independent of $x$. Since
$$
\int_0^L dx\,\rho^{(1)}_{N,L}(x)=N\,,
$$
denoting $N/L=\rho$, we get $\rho^{(1)}_{N,L}(x)=\rho$. Conditions (2.1)
imply $1/2\le\rho\le1$. Similarly
$\rho^{(2)}_{N,L}(x,y)$ depends only upon $|x-y|$, and
$$
\int_0^L dx\int_0^L dy\,\rho^{(2)}_{N,L}(x,y)=N(N-1)\,.
$$
The pair {\sl distribution} function is defined as
$$
g_{N,L}(x)={\rho^{(2)}_{N,L}(0,x)\over\rho^2}\eqno{(2.4)}
$$
and satisfies
$$
\int_0^L dx\,g_{N,L}(x)=\bigl(1-{1\over N}\bigr)\,L
$$
The degree of correlation or independence between the configuration around $x$
and the configuration around $y$ may be estimated by looking at
the pair {\sl correlation} function
$$
{\Pe_{N,L}(\,\exists i: X_i\in(x,x+dx),\exists j: X_j\in(y,y+dy)\,)
\over\Pe_{N,L}(\,\exists i: X_i\in(x,x+dx)\,)\,
\Pe_{N,L}(\,\exists j: X_j\in(y,x+dy)\,)}-1=g_{N,L}(x-y)-1
$$

\bigskip\noindent{\bf 3. Spacings}
\smallskip\noindent
Let $(X_1,X'_2,\dots,X'_N)$ be obtained by permutation of $(X_1,\dots,X_N)$ 
so that 
$$
X_1\le X'_2\le\dots\le X'_N\le X_1+L
$$ 
and let
$$
\eqalign{
S_1&=X'_2-X_1\cr
S_i&=X'_{i+1}-X'_i\,,\quad i=2,\dots N-1\cr
S_N&=X_1+L-X'_N
}
$$
Then $(X_1,S_1,\dots,S_{N-1})$ is distributed according to the Lebesgue
measure on the subset of $\Re^N$ defined by
$$
\eqalign{
0&\le X_1<L\cr
1&\le S_i\le2\,,\quad i=1,\dots N-1\cr
1&\le L-\sum_1^{N-1}S_i\le2
}\eqno{(3.1)}
$$        
and we may write (2.2) as
$$
\eqalign{
Z_{N,L}&=\int_0^L dx_1\int_1^2ds_1\dots\int_1^2ds_{N-1}\,
\un_{1\le L-\sum_1^{N-1}s_i\le2}\cr
&=L\int_1^2ds_1\dots\int_1^2ds_N\,\delta\Bigl(\,\sum_1^Ns_i-L\Bigr)
}\eqno{(3.2)}
$$
It follows that $(S_1,\dots,S_N)$ is distributed according to
$$
\eqalign{
\Pe_{N,L}\bigl(\,S_1\in(s_1,s_1+ds_1),\dots,S_N\in(s_N,s_N+ds_N)\,\bigr)=\cr
=\Bigl(\,\prod_{i=1}^N\un_{1\le s_i\le2}\,ds_i\,\Bigr)
\delta\Bigl(\,\sum_1^Ns_i-L\Bigr)\,\big/\,{\rm normalisation}
}\eqno{(3.3)}
$$
Then, for $x>0$,
$$
\eqalign{
\Pe_{N,L}(\,\exists i,j: X_i\in(0,dy),\,&X_j\in(x,x+dx)\,)
=N\Pe_{N,L}(\,X_1\in(0,dy), \exists j: X_j\in(x,x+dx)\,)\cr
&=N\sum_{m=1}^{N-1}\Pe_{N,L}(\,X_1\in(0,dy), X_1+S_1+\dots+S_m\in(x,x+dx)\,)\cr
&=N{dy\over L}\sum_{m=1}^{N-1}\Pe_{N,L}(\,S_1+\dots+S_m\in(x,x+dx)\,)
}
$$
so that
$$
g_{N,L}(x)=\sum_{m=1}^{N-1}{\Pe_{N,L}(\,S_1+\dots+S_m\in(x,x+dx)\,)\over\rho\, 
dx}\,,\qquad x>0\,,\eqno{(3.4)}
$$
and $g_{N,L}(-x)=g_{N,L}(x)$.

\bigskip\noindent{\bf 4. Isobaric ensemble}
\smallskip\noindent
In the thermodynamic limit, $L\to\infty$, $N\to\infty$ with $\rho=N/L$ fixed,
asymptotic statistical properties are easier to compute in an ensemble
where the pressure $p$ is fixed instead of the system size $L$. This is
defined as follows: for $p\in\Re$, let $\tilde S_1,\dots,\tilde S_N$ be
distributed according to
$$
\Pe_{N,p}\bigl(\,\ti S_1\in(s_1,s_1+ds_1),\dots,\ti S_N\in(s_N,s_N+ds_N)\,
\bigr)
=\Bigl(\,\prod_{i=1}^N\un_{1\le s_i\le2}\,\,e^{^{\scr -ps_i}}\,ds_i\,\Bigr)\,
\big/\,{\rm norm.}\eqno{(4.1)}
$$
The normalisation is
$$
Z_{N,p}=\Bigl(\,\int_1^2ds\,e^{^{\scr -ps}}\,\Bigr)^N
=\Bigl(\,{e^{^{\scr -p}}-e^{^{\scr -2p}}\over p}\,\Bigr)^N\eqno{(4.2)}
$$
and the associated potential is the Gibbs potential $G(N,p)=-\ln Z_{N,p}$.
The thermodynamic relation
$$
{\partial G\over\partial p}=\Ee_{N,p}\sum \ti S_i
$$
is the usual one, with $\Ee_{N,p}\sum \ti S_i=L$, the mean system size.
We thus get
$$
L={\partial G\over\partial p}={N\over p}
+N\,{e^{^{\scr -p}}-2e^{^{\scr -2p}}\over e^{^{\scr -p}}-e^{^{\scr -2p}}}
$$
or
$$
{1\over\rho}={1\over p}+{e^{^{\scr -p}}-2e^{^{\scr -2p}}\over e^{^{\scr -p}}
-e^{^{\scr -2p}}}\eqno{(4.3)}
$$
Although temperature plays only a dummy role, one may wish to have it in,
and also to replace the length {\sl one} of rods and maximum 
allowed gap between neighboring rods by a length $\ell$. Then one gets
$$
{1\over\rho}={kT\over p}+{e^{^{\scr -p\ell/kT}}-2e^{^{\scr -2p\ell/kT}}\over 
e^{^{\scr -p\ell/kT}}-e^{^{\scr -2p\ell/kT}}}\,,\eqno{(4.4)}
$$
which looks more like an equation of state.

Going back to (4.3), for $p=+\infty$, $\ln2$, $0$, $-\infty$, we have 
respectively $\rho=1,\,\ln2,\,2/3$, $1/2$.

In this ensemble we can compute, for $x>0$,
$$
\eqalign{
g_{N,p}(x)&=
\sum_{m=1}^{N-1}{\Pe_{N,p}(\,\ti S_1+\dots+\ti S_m\in(x,x+dx)\,)\over\rho\,dx}
\cr
&=\sum_{m=1}^{N-1}\biggl({p\over e^{^{\scr -p}}-e^{^{\scr -2p}}}\biggr)^m
\,{1\over\rho}\,\int_1^2e^{^{\scr -ps_1}}ds_1\dots\int_1^2e^{^{\scr -ps_m}}ds_m
\,\delta\Bigl(\,\sum s_i-x\,\Bigr)
\cr
&=\sum_{m=1}^{N-1}\biggl({p\over e^{^{\scr -p}}-e^{^{\scr -2p}}}\biggr)^m
\,{e^{^{\scr -px}}\over\rho}\int_0^1dt_1\dots\int_0^1dt_m
\,\delta\Bigl(\,\sum t_i-x+m\,\Bigr)\cr
&=\sum_{m=1}^{N-1}\biggl({p\over e^{^{\scr -p}}-e^{^{\scr -2p}}}\biggr)^m
\,{e^{^{\scr -px}}\over\rho}\,u_m(x-m)
}\eqno{(4.5)}
$$
where the densities $u_n(x)$ satisfy the recursion relation
$$
u_{n+1}(x)=\int_0^1u_n(x-y)\,dy\,,
$$
which gives [6, p. 27] 
$$
\eqalign{
u_1(x-1)&=\un_{x\in(1,2)}\cr
u_m(x-m)&={1\over (m-1)!}\sum_{\ell=0}^m(-)^\ell\biggl({m\atop \ell}\biggr)
\Bigl(x-m-\ell\Bigr)_+^{\ m-1}\,,\quad m\ge2
}\eqno{(4.6)}
$$
The thermodynamic limit $g(x)=\lim_{N\to\infty} g_{N,p}(x)$ is given by
(4.5)(4.6) with the range of $m$ extended to infinity. Note that for each $x$
only $[x/2]$ terms contribute, because $u_m(x-m)$ vanishes for $x>2m$.
Plots of $g(x)$ are given in Figure 1.

$\Pe_{N,p}(\cdot)$ also induces a distribution for parking configurations of
$N$ rods on $[0,\infty)$, so that $g_{N,p}(x)$ may be considered to be defined
by (2.3)(2.4) with $L$ replaced by $p$ and $1/\rho=\Ee_{N,p}L/N$. Then
$g_{N,p}(-x)=g_{N,p}(x)$ and $g(-x)=g(x)$.

$$
\psboxto(13cm;0cm){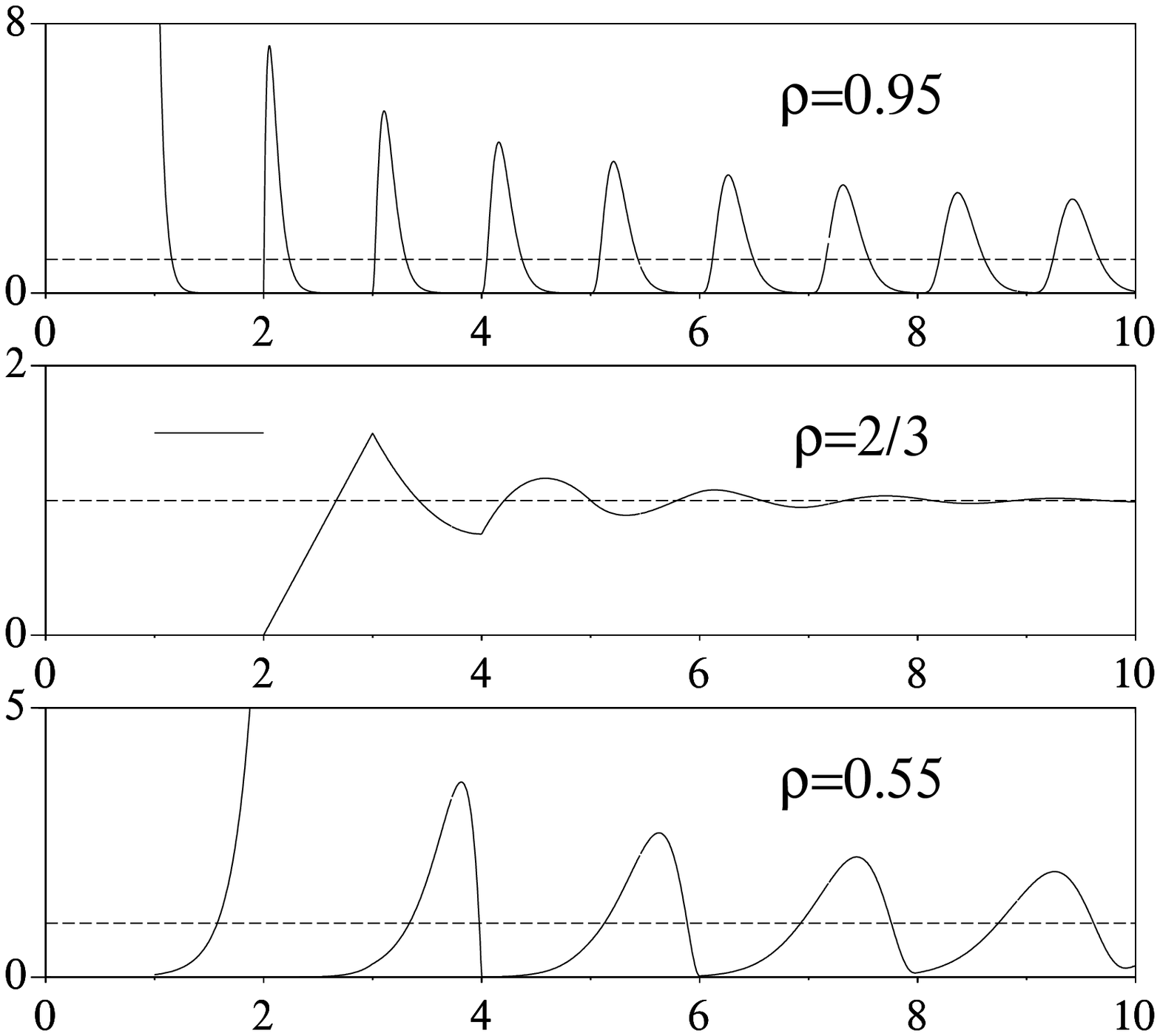}
$$
\centerline{\ninerm Fig. 1: The pair distribution function $g(x)$}

\bigskip\noindent{\bf 5. Correlation length}
\smallskip\noindent
Let us compute for $\la>0$ the Laplace transform,
$$
\eqalign{ 
\hat g^{\rm Laplace}(\la)&=\int_0^\infty e^{^{\scr -\la x}}g(x)\,dx\cr
&=
-{p\over\rho}\,{e^{^{\scr -(\la+p)}}-e^{^{\scr -2(\la+p))}}\over
-(\la+p)\bigl(e^{^{\scr -p}}-e^{^{\scr -2p}}\bigr)
+p\bigl(e^{^{\scr -(\la+p)}}-e^{^{\scr -2(\la+p)}}\bigr)}
}\eqno{(5.1)}
$$
Using (4.3), we find $\la\,\hat g^{\rm Laplace}(\la)\to 1$ as $\la\to 0$,
in agreement with $g(x)\to 1$ as $x\to\infty$.
We then consider the Laplace transform of the pair {\sl correlation} 
function $h(x)=g(x)-1$,
$$
\hat h^{\rm Laplace}(\la)=\hat g^{\rm Laplace}(\la)-{1\over\la}
$$
and look for its poles, which can only be in the complex half-plane
$Re(\la)\le0$. We must solve
$$
{e^{^{\scr -(\la+p)}}-e^{^{\scr -2(\la+p)}}\over\la+p}
={e^{^{\scr -p}}-e^{^{\scr -2p}}\over p}\eqno{(5.2)}
$$
Setting $\la=a+ib$ and $\ti p=p+a$, the modulus and the phase in (5.2) give
$$
\eqalign{
{\bigl(\,e^{^{\scr -\ti p}}-e^{^{\scr -2\ti p}}\,\bigr)^2
+2e^{^{\scr -3\ti p}}(1-\cos b)\over\ti p^2+b^2}
&={\bigl(\,e^{^{\scr -p}}-e^{^{\scr -2p}}\,\bigr)^2\over p^2}\cr
{\sin b-e^{^{\scr -\ti p}}\sin 2b\over
\cos b-e^{^{\scr -\ti p}}\cos 2b}
&=-{b\over\ti p}}\eqno{(5.3)}
$$
When $p\to+\infty$, looking only in the strip $-1\le a\le 0$, we get $b$ from 
the second equation and then $a$ from the first:
$$
\eqalign{
b&=2n\pi\bigl(\,1-{1\over p}+\O\bigl({1\over p^2}\bigr)\,\bigr)\cr
a&=-{2n^2\pi^2\over p^2}\,\bigl(\,1+\O\bigl({1\over p}\bigr)\,\bigr)
}\eqno{(5.4)}
$$
where $n\in\Ze\setminus\{0\}$.
When $p\to-\infty$, we get similarly
$$
\eqalign{
b&=n\pi\bigl(\,1-{1\over 2p}+\O\bigl({1\over p^2}\bigr)\,\bigr)\cr
a&=-{n^2\pi^2\over4p^2}\,\bigl(\,1+\O\bigl({1\over p}\bigr)\,\bigr)
}\eqno{(5.5)}
$$
The pair correlation function $h(x)$ is conveniently recovered by Fourier
transform: 
$$
\eqalign{
\hat h(k)=\int_{-\infty}^\infty e^{^{\scr ikx}}g(x)\,dx
=\hat h^{\rm Laplace}(ik)+\hat h^{\rm Laplace}(-ik)
}\eqno{(5.6)}
$$
$$
h(x)={1\over 2\pi}\int_{-\infty}^\infty e^{^{\scr -ikx}}\hat h(k)\,dk
$$
Using (5.1)(5.4)(5.5), the contour of integration can be shifted past
the nearest poles, whose residues are the dominant part of $h(x)$
when $x\to\infty$, giving the inverse correlation length 
$$
{1\over\xi}=\liminf_{x\to+\infty}-{1\over x}\ln |h(x)|=\left\{\matrix{
\dis{2\pi^2\over p^2}\,\bigl(\,1+\O\bigl({1\over p}\bigr)\bigr)
\qquad{\rm as}\quad p\to+\infty\cr
\dis{\pi^2\over 4p^2}\,\bigl(\,1+\O\bigl({1\over p}\bigr)\bigr)
\qquad{\rm as}\quad p\to-\infty}\right.\eqno{(5.7)}
$$
or, using (4.3),
$$
\xi\sim\left\{\matrix{
\dis{1\over2\pi^2}\,(1-\rho)^{-2}\qquad{\rm as}\quad &\rho\nearrow1\cr\cr
\dis{1\over\pi^2}\,(2\rho-1)^{-2}\qquad{\rm as}\quad &\rho\searrow1/2
}\right.\eqno{(5.8)}
$$
The imaginary part of the poles gives the pseudo-period of oscillation,
$\Delta x=1$ when $\rho\nearrow1$ and $\Delta x=2$ when $\rho\searrow1/2$,
in accordance with Figure 1 and with the Dirac train limits
$$
\lim_{p\to\pm\infty}\hat g^{\rm Laplace}(\la)=\left\{\matrix{
{e^{^{\scr -\la}}\over1-e^{^{\scr -\la}}}
\qquad{\rm as}\quad &p\to+\infty\cr\cr
{2e^{^{\scr -2\la}}\over1-e^{^{\scr -2\la}}}
\qquad{\rm as}\quad &p\to-\infty
}\right.\eqno{(5.9)}
$$
$$
\lim_{p\to\pm\infty} g(x)=\left\{\matrix{\dis
\sum_{\dis n\in \Ze}\delta(x-n)
\qquad{\rm as}\quad &p\to+\infty\cr\cr
2\dis\sum_{\text n\in \Ze}\delta(x-2n)
\qquad{\rm as}\quad &p\to-\infty
}\right.\eqno{(5.10)}
$$

\bigskip\noindent{\bf 6. Equivalence of ensembles}
\medskip\noindent
For any $p\in\Re$ we may write (3.2) as
$$
\eqalign{
{Z_{N,L}\over L}
&=\int_1^2ds_1\dots\int_1^2ds_N\,\delta\Bigl(\,\sum s_i-L\Bigr)\cr
&=\int_1^2ds_1\dots\int_1^2ds_N\,
e^{^{\scr -p\,\bigl(\sum s_i-L\bigr)}}\,\delta\Bigl(\,\sum s_i-L\Bigr)\cr
&=e^{^{\scr pL}}\,Z_{N,p}\ \Ee_{N,p}\ \delta\Bigl(\,\sum\ti S_i-L\Bigr)
}\eqno{(6.1)}
$$
where $\Ee_{N,p}(\cdot)$ is defined by (4.1)(4.2).  Let us choose $p$ so that
$$
\Ee_{N,p}\sum\ti S_i=L
$$
where $L$ is the fixed value in the fixed $L$ ensemble.
Then the Central Limit Theorem implies
$$
\Ee_{N,p}\ \delta\Bigl(\,\sum\ti S_i-L\Bigr)=\O(L^{-1/2})
$$
and we get
$$
{Z_{N,L}\over L}=e^{^{\scr pL}}\,Z_{N,p}\ \O(L^{-1/2})
$$
The free energy per particle at density $\rho$,
$$
f(\rho)=\lim_{N\to\infty}{F(N,{N/\rho})\over N}
$$
and the Gibbs potential per particle at the corresponding pressure $p$,
$$
g(p)=\lim_{N\to\infty}{G(N,p)\over N}=
-\ln{e^{^{\scr -p}}-e^{^{\scr -2p}}\over p}\eqno{(6.2)}
$$
are therefore related by
$$
f(\rho)=-{p\over\rho}+g(p)\eqno{(6.3)}
$$
The local version of the central limit theorem can then be used to prove
$$
\lim_{N\to\infty} g_{N,N/\rho}(x)=\lim_{N\to\infty} g_{N,p}(x)=g(x)
$$
along the lines of the proof of the Wulff shape for one-dimensional
interfaces in [5]. Equivalence of ensembles is of course very standard,
but its justification or derivation is generally more involved than for
the present model as outlined above. 

The explicit form of the free energy $f(\rho)$ allows to answer an old question
in statistics [7,8]: let $X_1,\dots,X_N$ be $N$ independent random variables
each distributed uniformly over the interval $(0,L)$. Denote
$\Pe_{N,L}^{\rm\, free}$ the corresponding probability distribution (the ideal 
gas), and call {\sl Parking} the event that the smallest and largest spacings
are respectively larger than one and smaller than two. We have
$$
\Pe_{N,L}^{\rm\, free}({\sl Parking})
={\int_1^2ds_1\dots\int_1^2ds_N\ \delta\bigl(\,\sum s_i-L\bigr)
\over\int_0^Lds_1\dots\int_0^Lds_N\ \delta\bigl(\,\sum s_i-L\bigr)}
={Z_{N,L}/L\over L^{N-1}/(N-1)!}\eqno{(6.4)}
$$
Let $1/2<\rho<1$. Then (6.1)(6.2)(6.3)(6.4) and Stirling's formula give 
$$
\lim_{N\to\infty}-{1\over N}\ln\Pe_{N,N/\rho}^{\rm\, free}({\sl Parking})
=1-{p\over\rho}-\ln\rho-\ln{e^{^{\scr -p}}-e^{^{\scr -2p}}\over p}\eqno{(6.5)}
$$
where $p$ is related to $\rho$ by (4.3).
\bigskip
\noindent{\bf References}
\smallskip
\item{[1]} A. R\'{e}nyi,  {\sl On a one-dimensional problem
concerning random space filling}, Magyar Tud. Akad. Mat. Kutato Int. Kozl,
{\bf 3}, 109-127 (1958); Selected papers of Alfr\'ed R\'{e}nyi, Vol. 2,
pp. 173-188, Akad\'emiai Kiad\'o, Budapest (1976).

\item{[2]} J. Frenkel, {\sl Kinetic Theory of Liquids}, Oxford University 
Press, New-York (1940), p.126.

\item{[3]} F. G\"ursey, {\sl Classical statistical mechanics of a rectilinear 
assembly}, Proc. Cambridge Phil. Soc. {\bf 46}, 182--194 (1950).

\item{[4]} Z. Salsburg, R. Zwanzig, J. Kirkwood, {\sl Molecular Distribution
Functions in a One-Dimensional Fluid}, J. Chem. Phys. {\bf 21}, 1098-1107
(1953).

\item{[5]} J. De Coninck, F. Dunlop, V. Rivasseau, {\sl On the
Microscopic Validity of the Wulff Construction and of the Generalized
Young Equation}, Commun. Math. Phys. {\bf 121}, 401--419 (1989). 

\item{[6]} W. Feller, {\sl An introduction to probability
theory and its applications}, {Vol. 2}, John Wiley and Sons, Second Edition,
New York (1971).

\item{[7]} D. A. Darling, {\sl On a class of problems related to the random 
division of an interval}, Ann. Math. Statistics, {\bf 24}, 239--253 (1953).

\item{[8]} Paul L\'evy, {\sl Sur la division d'un segment par des points 
choisis au hasard}, C. R. Acad. Sci. Paris, {\bf 208}, 147--149 (1939).

\bye